\documentclass[12pt]{article}
\usepackage{amsmath}
\usepackage{amssymb}
\usepackage{amsfonts}
\usepackage{eucal}
\usepackage[usenames]{color}
\usepackage{graphicx}
\numberwithin{equation}{section}
\newcommand{\R}{\mathbb{R}}
\newcommand{\E}{\mathbb{E}}

\newtheorem{Theorem}{Theorem}[section]

\newtheorem{Remark}[Theorem]{Remark}
\begin{document}
\title{On the first-passage area of a L$\acute{\text{e}}$vy process}
\author{Mario Abundo$^*$, Sara Furia \thanks{Dipartimento di Matematica, Universit\`a  ``Tor Vergata'', via della Ricerca Scientifica, I-00133 Rome, Italy.
Corresponding author E-mail: \tt{abundo@mat.uniroma2.it}}
}
\date{}
\maketitle

\begin{abstract}
\noindent Let be $X(t)= x - \mu t + \sigma B_t - N_t$ a L$\acute{\text{e}}$vy process starting from $x >0,$
where $ \mu \ge 0, \ \sigma \ge 0,  \ B_t$ is a standard BM, and $N_t$ is a homogeneous Poisson process with intensity
$ \theta >0,$ starting from zero.
We study the joint distribution of the first-passage time below zero, $\tau (x),$ and the first-passage
area, $A(x),$ swept out by $X$ till the time $\tau (x).$  In particular, we establish differential-difference
equations with outer conditions for the Laplace transforms of $\tau(x)$ and $A(x),$ and for their joint moments.
In a special case $(\mu = \sigma =0),$ we show
an algorithm to find recursively the  moments $E[\tau(x)^m A(x)^n],$ for any integers $m$ and $n;$ moreover, we obtain the expected value of the
time average of $X$ till the time $\tau(x).$
\end{abstract}


\section{Introduction}
This is a continuation of the articles \cite{abundo:mcap13} and \cite{abdelve:mcap17}; actually, in \cite{abundo:mcap13}
we studied the distribution
of the first-passage area (FPA)  $A(x)= \int _0 ^ { \tau (x)} X(t) dt ,$ swept out by a one-dimensional jump-diffusion process $X(t),$ starting from $x>0,$
till its first-passage time (FPT) $\tau(x)$
below zero, while in \cite{abdelve:mcap17} we examined the special case (without jumps) when
$X(t)$ is Brownian motion (BM) $B_t$  with negative drift $- \mu,$ that is, $X(t)= x - \mu t + B_t, $
studying in particular the joint distribution of $\tau(x)$ and $A(x).$
In the present paper, we aim to investigate the  joint distributions of $\tau(x)$ and  $A(x),$ in the case when
$X(t)$ is a L$\acute{\text{e}}$vy process of the form
\begin{equation} \label{levyprocess}
X(t)= x - \mu t + \sigma B_t - N_t , \ x >0 ,
\end{equation}
where $ \mu \ge 0, \ \sigma \ge 0,  \ B_t$ is standard BM, and $N_t$ is a homogeneous Poisson process with intensity
$ \theta >0,$ starting from zero; thus, $X(t)$ turns out to be the superposition  of drifted BM
and Poisson process. \par\noindent
Referring to the L$\acute{\text{e}}$vy process \eqref{levyprocess},
we state and solve differential-difference equations for the Laplace transform
of the two-dimensional random variable $(\tau(x), A(x)).$ In particular, when $\mu = \sigma =0,$ we obtain
the joint moments  $E[\tau(x)^m A(x)^n]$ of the FPT and FPA, and we present an  algorithm  to find them recursively,
for any $m$ and $n;$
moreover, we find the expected value of the time average of $X(t)$ till its FPT below zero. \par
Studying the FPA of a process such as  \eqref{levyprocess} is peculiar when modeling the evolution of certain random systems described by the superposition of a continuous stochastic process and a jump process (see references in \cite{abundo:mcap13}); these arise e.g.
in solar physics studies, non-oriented  animal movement patterns, and
DNA breathing dynamics, as regards systems where the jump component can be absent (see e.g. \cite{kearney:jph07}, \cite{knight:jam00},
\cite{perman:aap96} and
references in \cite{kearney:jph14}). Applications can be found
in {\it Queueing Theory} (see e.g. \cite{abundo:mcap13}),
%if one identifies $X(t)$ with the length, at time $t,$ of a queue  subjected to possible sudden changes, and $\tau(x)$ with the busy period, that is the time %until the queue is first empty; then, $A(x)$ represents the cumulative waiting time experienced by all the ``customers''  during a busy period. Another
and in {\it Finance}, in the framework of default-at-maturity model, which assumes the exchange rate follows a jump-diffusion process
(see e.g. \cite{bates:1996}, \cite{jorion:1988}); for other examples from Economics and Biology, see e.g. \cite{abundo:mcap13}.
%if the variable $t$ represents the quantity of a commodity that producers have available for sale, and $X(t)$ describes the price of the commodity as a function %of the quantity in a supply-and-demand model; then, one can identify $\tau(x)$ with the amount of product at which the price falls to practically zero, while the %area $A(x)$ under the demand curve provides a measure of the total value that consumers receive from consuming that amount of the product.
\par
The paper is organized as follows. Section 2 contains preliminary
results on jump-diffusion processes. In section 3 we deal with the
Laplace transform of $(\tau(x), A(x))$ and the joint moments of
$\tau(x)$ and $A(x),$ in the case when $X(t)$ is a Poisson
process (that is $\mu = \sigma =0)$; precisely, we
find explicitly $E[\tau(x) A(x)],$ and we establish ODEs for the
joint moments
 $E[\tau(x)^m A(x)^n]$ of order $n+m, \ n, m \ge 0,$  presenting also
an algorithm to find recursively them; moreover, we
find the expected value of the time average of $X(t)$ till its FPT below zero.
Sections 4 and 5 are respectively devoted to study the distributions of $\tau(x), \ A(x)$ in the case of Poisson process with negative drift ($\sigma =0),$  and drifted BM with Poisson jumps ($\mu \neq 0, \ \sigma \neq 0).$ Finally, Section 6 contains conclusions and final remarks.

\section{Notations, formulation of the problem and preliminary results}
We recall from \cite{abundo:mcap13} some definitions and results concerning FPT and FPA for a \par\noindent
one-dimensional, time homogeneous jump-diffusion process $X(t)$ driven by the SDE:
\begin{equation} \label{jumpdiffueq}
 dX(t) = b(X(t))dt +\sigma(X(t)) dB_t + \int _{-\infty}^{+\infty} \gamma
(X(t),u) \nu(dt,du)
\end{equation}
with assigned initial condition $X(0) = x >0 ;$ here $\nu (\cdot
,\cdot)$ is a temporally homogeneous Poisson random measure (see
\cite{gimsko:sde72} for the definitions), and the functions
$b(\cdot), \sigma (\cdot), \gamma (\cdot, \cdot)$ satisfy suitable
conditions for the existence and  uniqueness of the solution (see
\cite{gimsko:sde72}, \cite{ikwa:sde81}). The random measure $\nu$
is supposed to be homogeneous with respect to time translation,
that is, its intensity measure $E(\nu(dt, du))$ is of the form $ E
[ \nu  (dt,du  ) ] = dt  \pi (du) $ for some positive measure
$\pi$ defined on $ {\cal B} (\rm I\!R),$
the Borel $\sigma-$field of subsets of $\rm I\!R,$ and
the jump intensity $
 \Theta = \int _{-\infty}^ {+\infty} \pi (du)  \ge 0
$
is assumed to be finite.
\par\noindent
If $\gamma = 0,$  or $\nu = 0,$
then the SDE \eqref{jumpdiffueq} becomes the usual It${\rm \hat{o}}$'s stochastic
differential equation for a simple-diffusion (i.e. without jumps).
If, for instance, the  measure $\pi$ is concentrated
over the set $\{u_1,u_2\}=\{-1,1\}$ with $\pi(u_i)=\theta _i$
and $\gamma(u_i)=\epsilon _i,$ the SDE
\eqref{jumpdiffueq} assumes the form
$ dX(t)=b(X(t))dt + \sigma(X(t)) dB_t + \epsilon _2 dN_2 (t) + \epsilon _1
dN_1 (t), $
where $ \epsilon _1 <0, \ \epsilon _2 >0 $ and $N_ i(t), \  t \ge 0 $ are independent
homogeneous Poisson processes of
amplitude $1 $ and rates $\theta _1$ and
$\theta _2,$ respectively  governing downward $(N_1)$ and upward $(N_2)$
jumps. \par
Let $D$ be the class of function $f(x) \in C^2,$
for which the function
$f(x+\gamma(x,u)) - f(x)$ is $\pi-$integrable for any $x.$
The differential operator associated to the process $X(t)$
is defined for  $f \in D$ by:
\begin{equation} \label{generator}
L f (x) = \frac 1 2
\sigma ^2 (x) \frac {d ^2 }  {d x^2} f(x)+
b(x) \frac { d }  {d x} f(x) + \int _ {- \infty} ^{+ \infty} [f(x+ \gamma(x,u))-f(x)] \pi(du) .
\end{equation}
%Then, from  the generalized It$\hat {\rm o}$'s formula  (see \cite{gimsko:sde72}), taking expectation, one obtains
%\begin{equation}
%E[f(t,X(t))] = f(0,X(0)) + E \left ( \int _0 ^t \left [\frac {\partial f}
%{\partial s} (s,X(s)) + L f (s,X(s)) \right ] ds \right ).
%\end{equation}
Let us define, for $x>0:$
\begin{equation} \label{firstpassagetime}
\tau(x) = \inf \{ t >0: X(t) \le 0 | X(0)=x \},
\end{equation}
that is the first-passage time below zero of $X(t),$ and assume that $\tau(x)$ is finite with
probability one. Really, it is possible to show (see \cite{abundo:pms00},  \cite{tuckwell:jap76}) that the probability $p_0(x)$ that
$X(t)$ ever leaves the interval $(0, + \infty)$ satisfies the
partial differential-difference equation (PDDE)
$
L p_0 = 0
$
with  outer condition:
$ p_0(y) =1 \ {\rm if}  \  y \le 0 .$
The equality $p_0(x) \equiv 1$ is equivalent to say that $\tau(x)$ is finite with probability one.
For diffusion processes without jumps (i.e. $\gamma =0)$
sufficient conditions are also available which ensure that $\tau(x)$ is
finite w.p. 1, and they concern the convergence of certain
integral associated to the coefficients of \eqref{jumpdiffueq}
(see e.g. \cite{gimsko:sde72}, \cite{has:sto80}).
\par
Let $U$ be a functional of the process $X;$ for $\lambda >0$ denote by
\begin{equation}
M_ {U, \lambda} (x) =  E \left [ e^ { - \lambda
\int _0^ {\tau(x)} U(X(s)) ds} \right ]
\end{equation}
the Laplace transform of the integral $\int _0 ^ {\tau(x)}
U(X(s))ds.$ Then, it holds (see \cite{abundo:mcap13}):
\begin{Theorem} \label{laplacetransform}
Let $X(t)$ be the solution of the SDE \eqref{jumpdiffueq}, starting from $X(0)=x >0;$ then, under the above assumptions, $M_{ U, \lambda} (x)$ is the solution of
the problem with outer conditions:
\begin{equation} \label{problemlaplace}
\begin{cases}
L M_{ U, \lambda} (x) = \lambda U(x) M_{U, \lambda } (x ) \\
M_ {U, \lambda }(y) =1,  \ {\rm for} \ y \le 0,  \
\lim _ { x \rightarrow + \infty} M_{U, \lambda } (x) =0,
\end{cases}
\end{equation}
where $L$ is the generator of $X,$ which is defined by \eqref{generator}.
\end{Theorem}
\hfill $\Box$
\bigskip
\par\noindent
The n-th order moment of $\int _0 ^ {\tau(x)}
U(X(s))ds,$ if it exists finite, is given by $(n=1, 2, \dots ):$
$$ T_ n (x) = E \left [  \left ( \int _0 ^ {\tau(x)}
U(X(s))ds \right ) ^n \right ]  = (-1)^n \left [ \frac {\partial
^n } {\partial \lambda ^n } M_ {U, \lambda } (x) \right ] _ { \lambda =0} .$$
Then, taking the n-th derivative with respect to $\lambda$ in
both members of the  equation  \eqref{problemlaplace}, and
calculating it for $\lambda =0,$ one easily obtains that
%\begin{Proposition} \label{propositionmoments}
the n-th order moment $T_n(x) \ (n=1, 2, \dots )$ of $\int _0 ^ {\tau(x)} U(X(s)) ds, $ whenever it exists finite,  is the solution of the PDDE:
\begin{equation} \label{eqmoments}
LT_n (x) = -n U(x) T_{n-1} (x) , \ \  x >0,
\end{equation}
which satisfies
\begin{equation} \label{outerconditionmoments}
 T_n(x)= 0, \ {\rm for } \ x \le 0
\end{equation}
and an  appropriate  additional  condition. \par
\noindent Indeed, the only condition $T_n (x)=0$ for $x \le 0$ is not sufficient to determinate uniquely the desired solution of
the PDDE \eqref{eqmoments}, because it is a second
order equation.  Note that for a diffusion without jumps $(\gamma =0 )$ and
for $U(x)\equiv 1,$
 \eqref{eqmoments} is nothing but the celebrated Darling and Siegert's equation  (\cite{darling:ams53}) for the moments of the first-passage time,
 and \eqref{outerconditionmoments}
 becomes simply the boundary condition $T_n(0) =0.$
\bigskip

In the next section, we start considering
the Poisson process $X(t)=x - N_t  \ (x>0),$ obtained from \eqref{levyprocess} by taking $\mu = \sigma =0,$
which is a special case of
the jump-diffusion driven by \eqref{jumpdiffueq}. As easily seen, the FPT of $X(t)$ below zero, $\tau(x),$  is
finite with probability one.
The successive sections regard Poisson process with negative
drift $X(t)= x - \mu t - N_t \ (\sigma =0)$ and drifted BM with Poisson jumps $X(t)=  x - \mu t + \sigma B_t - N_t , \ (\sigma >0).$
\par\noindent
To illustrate the behavior of the process $X(t)$ in \eqref{jumpdiffueq}, in the figure below we show a sample path of the process
$X(t),$ respectively in the case of Poisson process, Poisson process with drift,
and drifted BM with Poisson jumps.
\begin{figure}[h]
\centering
{\includegraphics[width=0.65\textwidth]{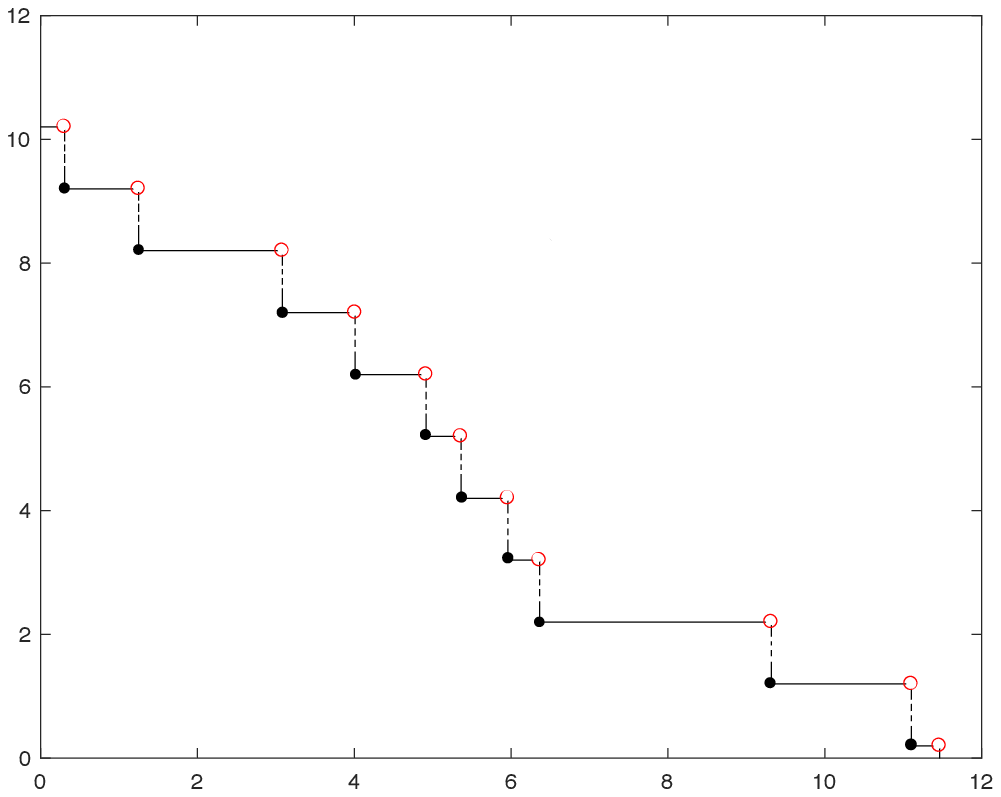}}
{\includegraphics[width=0.65\textwidth]{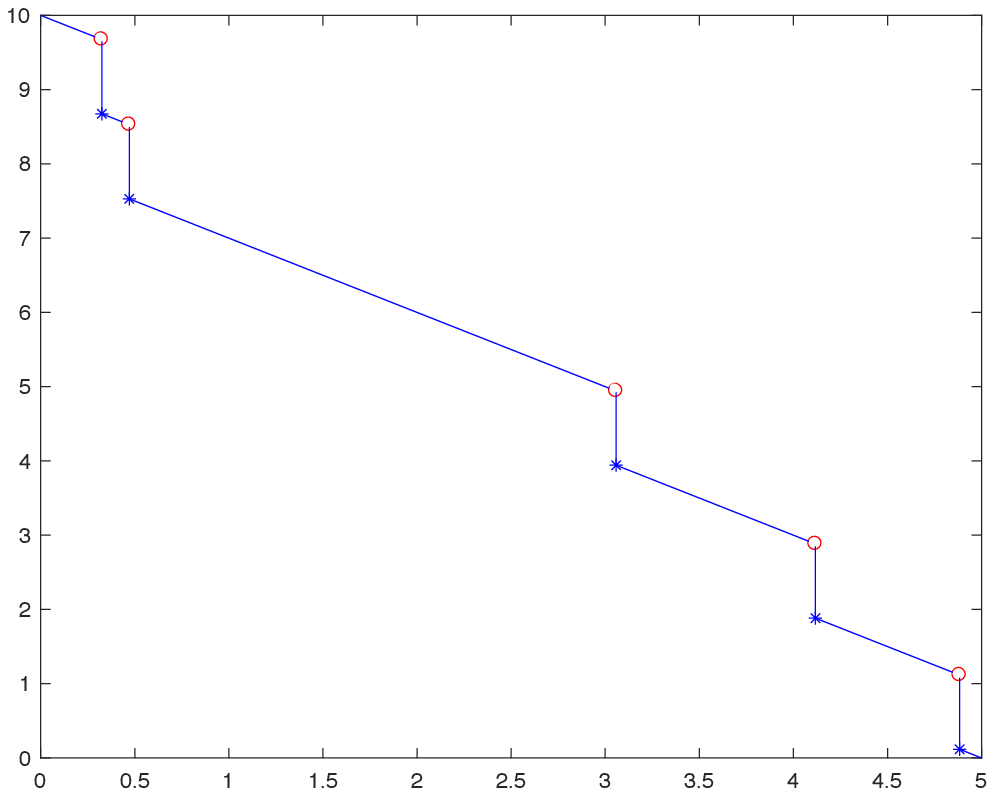}}
{\includegraphics[width=0.65\textwidth]{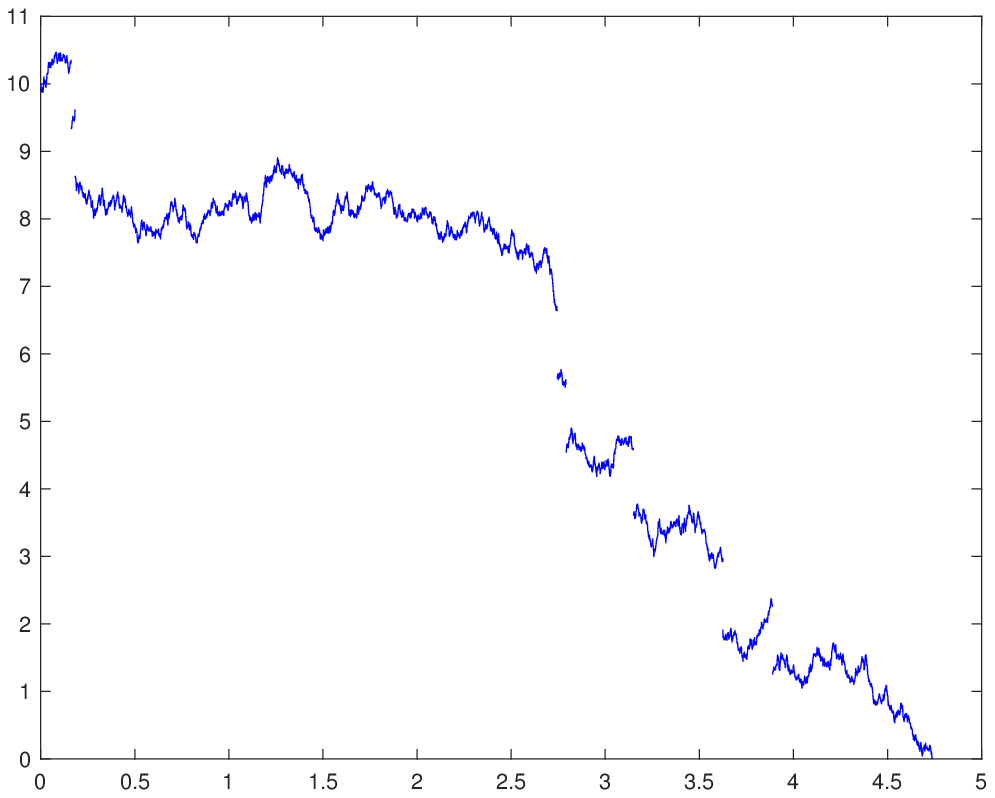}}
\caption{A sample path of the process $X(t) = x - N(t)$ with $x=10.2$ and $\theta = 1 $(top); $X(t)=x - \mu t - N_t $
with  $x=10, \ \theta = 1$ and $ \mu =1 $
(middle); $X(t)= x - \mu t + \sigma  B_t - N_t ,$ with $x=10, \ \theta = 1$, $ \mu =1 $ and $ \sigma =1 $ (bottom).}
\label{paths}
\end{figure}

\section{The case of Poisson process $(\mu = \sigma =0)$}
For $x>0,$ let us consider the process $X(t)= x - N_t,$ where $N_t \ (N_0=0) $ is a homogeneous Poisson process with intensity
$\theta >0.$
The infinitesimal generator is:
\begin{equation}
L f(x) = \theta [f(x-1) - f(x)] , \ f \in C^0(\rm I\!R),
\end{equation}
and
$ \tau (x) = \inf \{ t >0 : x- N_t \le  0 \}, \ A(x)= \int _ 0 ^{ \tau (x)} (x - N_t ) dt .$ \par
By Theorem \ref{laplacetransform} with $U(x)=1,$ it follows that the
Laplace transform $M_ {U, \lambda } (x)$ of $\tau (x)$ is the solution of the equation $L M _ {U,  \lambda } (x) = \lambda M _ {U,  \lambda } (x),$ with
outer condition $M_ {U, \lambda } ( y) =1$ for $ y \le 0.$ By solving this equation, we get for any $x >0$ (see \cite{abundo:mcap13}):
\begin{equation} \label{LTtaupoisson}
 M_ {U, \lambda } (x) =
\begin{cases}
\left ( \frac {\theta } {\theta + \lambda } \right ) ^ x  &   {\rm if} \ x \in \rm I\!N    \\
\left ( \frac {\theta } {\theta + \lambda } \right ) ^ {[x]+1}  &   {\rm if} \ x \notin \rm I\!N ,
\end{cases}
\end{equation}
where $[x]$ denotes the integer part of $x.$
Recalling the expression of the Laplace transform of the Gamma density, we note that $\tau(x)$ has Gamma distribution
with parameters $(x, \theta)$ if $x$ is a positive integer, while it has Gamma distribution
with parameters $([x]+1, \theta)$ if $x$ is not an integer.
In Fig. 2, we show the graph of the Laplace transform of $\tau(x),$ as a function of $ \lambda >0,$ for some values of the parameters. \par
\begin{figure}[ht]
\centering
{\includegraphics[width=0.55 \textwidth]{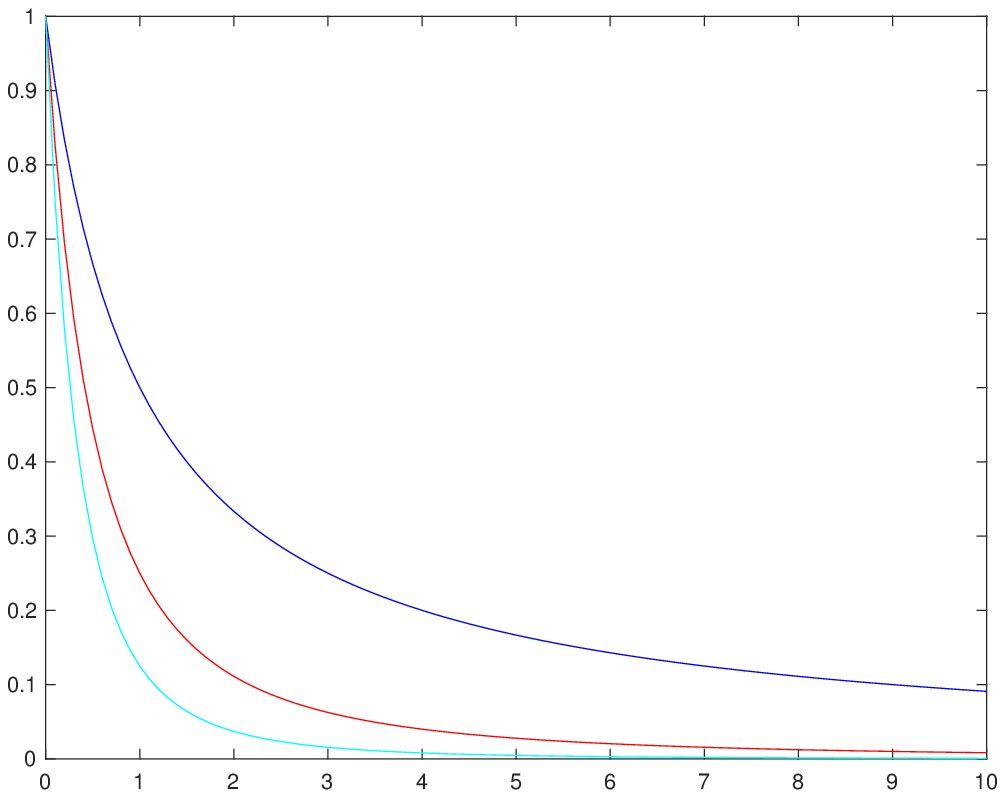}}{\includegraphics[width=0.55\textwidth]{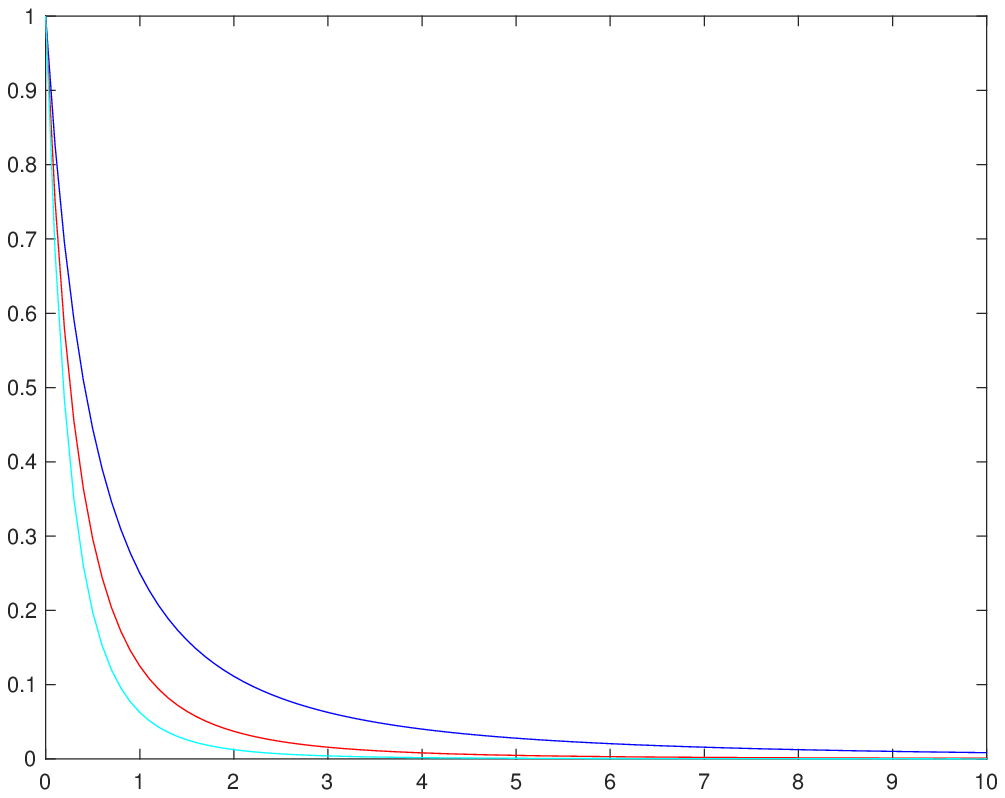}}
\caption{Laplace Transform of $\tau(x),$ as a function of $\lambda >0 ,$ for $\theta =1;$  on the left, from top to bottom:
$x = 1$ (blue), $x = 2$ (red), and $x = 3$ (cyan);
on the right, from top to bottom: $x = 1.5$ (blue), $x = 2.5$ (red) and $x = 3.5$ (cyan).}
\end{figure}
The moments
$ T_ n (x) = E \left [ ( \tau (x))^n \right ]$ are soon obtained by the formula  \par\noindent
$T_n (x) =   (-1)^n \left [ \frac {\partial
^n } {\partial \lambda ^n } M_\lambda (x) \right ] _ { \lambda =0}.$
In fact, we get (see \cite{abundo:mcap13}):
\begin{equation} \label{primiduemomtaupoisson}
E ( \tau(x))=
\begin{cases}
\frac {x } {\theta } &  {\rm if} \ x \in \rm I\!N    \\
\frac {[x]+1 } {\theta}  &  {\rm if} \ x \notin \rm I\!N
\end{cases}, \ {\rm and } \
E ( \tau ^2(x))=
\begin{cases}
\frac {x^2 } {\theta ^2 } + \frac {x } {\theta ^2 }   &  {\rm if} \ x \in \rm I\!N    \\
\frac {([x]+1)^2 } {\theta ^ 2} + \frac {[x]+1 } {\theta ^2  }   &  {\rm if} \ x \notin \rm I\!N .
\end{cases}
\end{equation}
Therefore:
$$ Var( \tau(x))=
\begin{cases}
\frac {x } {\theta ^2 }  &  {\rm if} \ x \in \rm I\!N    \\
\frac {[x]+1 } {\theta ^2 }   &  {\rm if} \ x \notin \rm I\!N .
\end{cases}
$$
\indent By Theorem \ref{laplacetransform} with $U(x)=x,$ we get the
Laplace transform $M_ {U, \lambda } (x)$ of $A(x)$ as the solution of the equation $L M _ {U,  \lambda } (x) = \lambda x M _ {U,  \lambda } (x),$ with
outer condition $M_ {U, \lambda } ( y) =1$ for $ y \le 0.$ By solving this equation, we obtain for any $x >0$
(see \cite{abundo:mcap13}):
\begin{equation} \label{LTApoisson}
M_ {U, \lambda } (x) =
\begin{cases}
\theta ^x \cdot \{(\theta + \lambda) ( \theta + 2 \lambda) \cdots (\theta + x \lambda)\}  ^ {-1}    &   {\rm if} \ x \in \rm I\!N    \\
\theta ^ {[x]+1} \cdot \{(\theta + \lambda x) ( \theta + \lambda (x-1)) \cdots (\theta + \lambda (x - [x])) \}^ {-1} &   {\rm if} \ x \notin \rm I\!N
\end{cases}
\end{equation}
that can be written in the unique form, valid for any $x >0:$
$$ M_ {U, \lambda } (x) = \frac { \theta ^{[x]+1}} {(\theta + \lambda x) ( \theta + \lambda (x-1)) \cdots
(\theta + \lambda (x - [x])) }  \ .$$
Notice that $M_ {U, \lambda } (x)$ turns out
to be the Laplace transform of a linear combination of $[x]+1$
independent exponential random variables with parameter $\theta,$
with coefficients $x, x-1, \dots , x-[x] .$  In the Figure 3, we show the graph of the Laplace transform of $A(x),$ as a function of
$ \lambda >0,$ for some values of the parameters.

\begin{figure}[ht]
\centering
{\includegraphics[width=0.55 \textwidth]{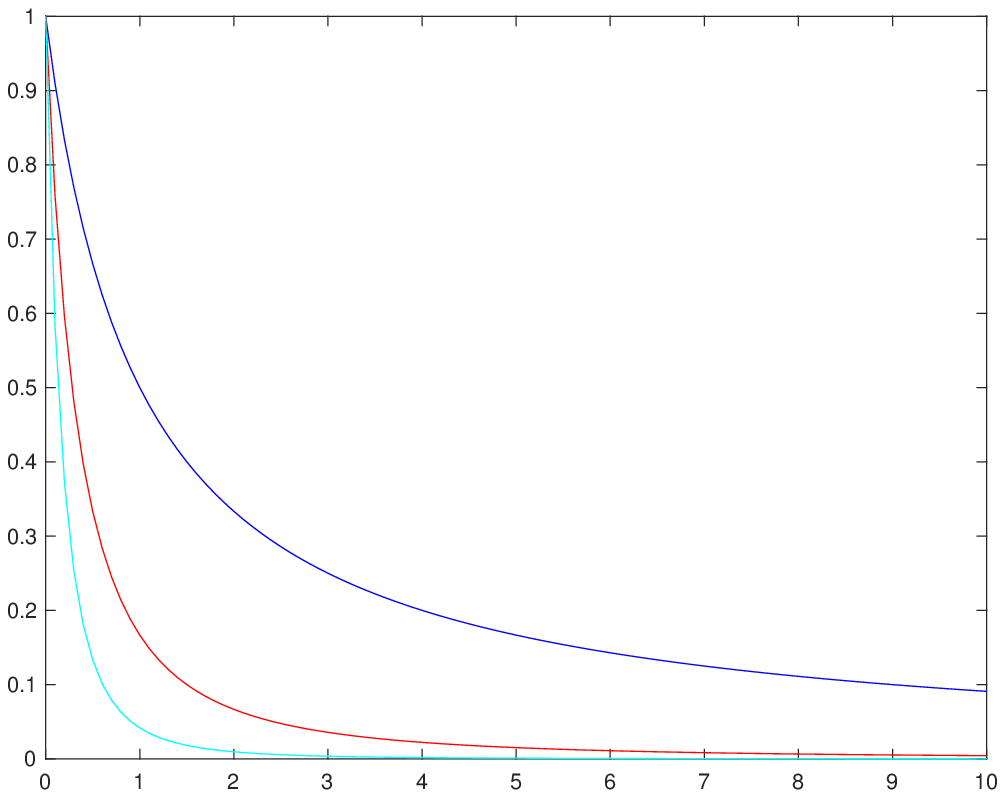}}{\includegraphics[width=0.55 \textwidth]{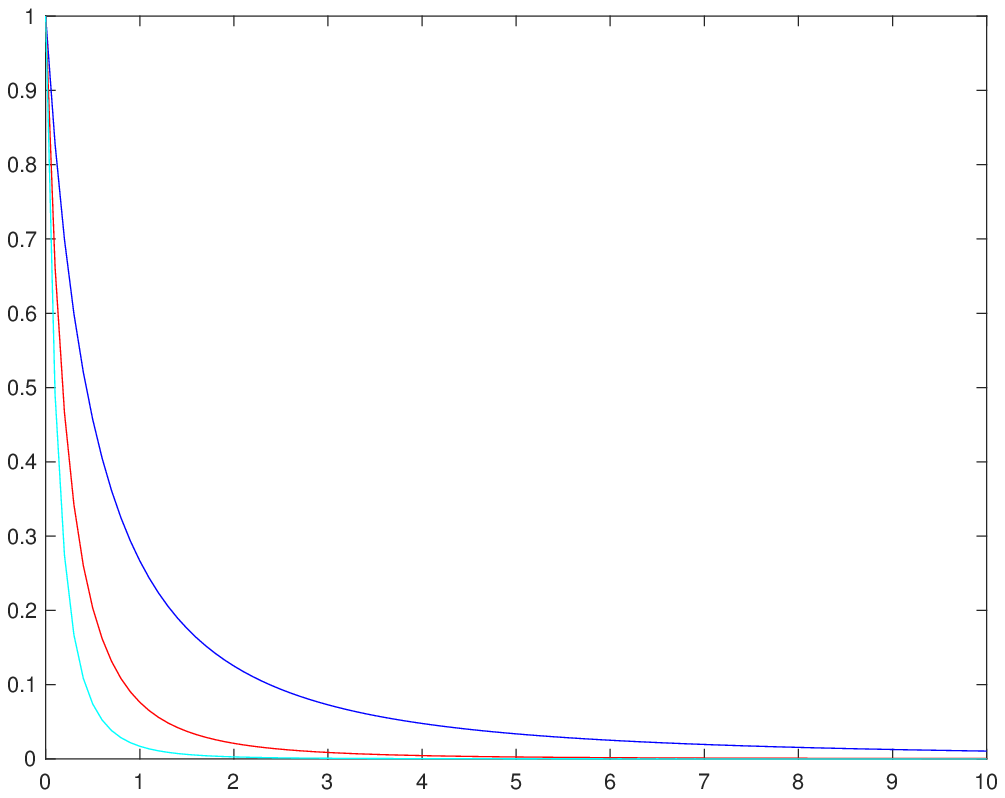}}
\caption{Laplace Transform of $A(x),$ as a function of $\lambda >0,$ for $\theta =1;$
on the left, from top to bottom: $x = 1$ (blue), $x = 2$ (red), and $x = 3$ (cyan);
on the right, from top to bottom: $x = 1.5$ (blue), $x = 2.5$ (red), and $x = 3.5$ (cyan).}
\end{figure}

The $n-$th order
moment  of $A(x)$ is given by $(-1)^n \left [ \frac {\partial ^n }
{\partial \lambda ^n } M_ {U, \lambda } (x) \right ] _ { \lambda
=0};$ calculating the first and second  derivative, we obtain (see \cite{abundo:mcap13}):
\begin{equation} \label{mediaApoisson}
 E(A(x)) = \frac
{(2x -[x])([x]+1) } { 2 \theta}
\end{equation}
 and
\begin{equation} \label{secondomomApoisson}
E ( A^2(x))= \frac {[x]+1 } {12 \theta ^2 } \Big ( 12x(x -[x])([x]+2) + [x] (3[x]^2+7 [x] +2)\Big )  .
\end{equation}

In the Figure 4, we show the densities of $\tau (x)$ and $A(x),$ estimated by Monte Carlo simulation, for some values of
the parameters.

\begin{figure}[ht]
\centering
{\includegraphics[width=0.6 \textwidth]{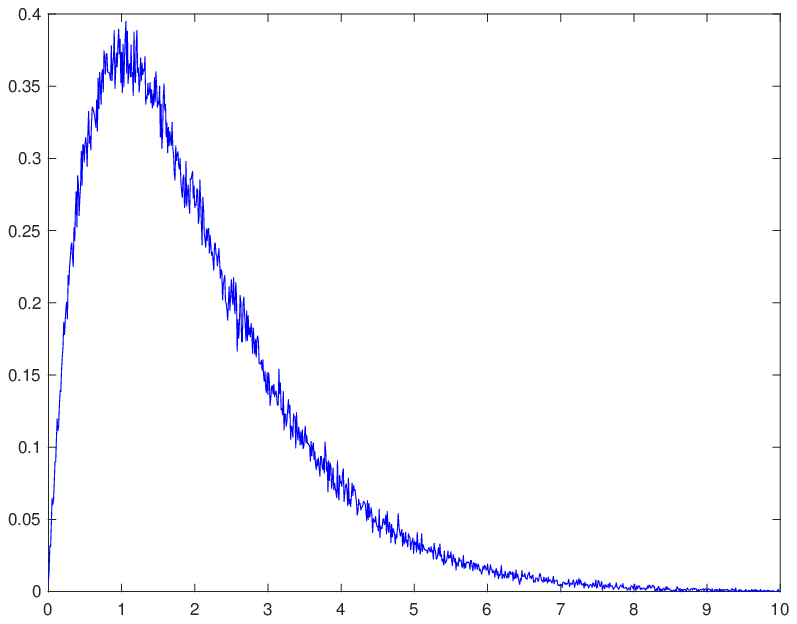}}{\includegraphics[width=0.6 \textwidth]{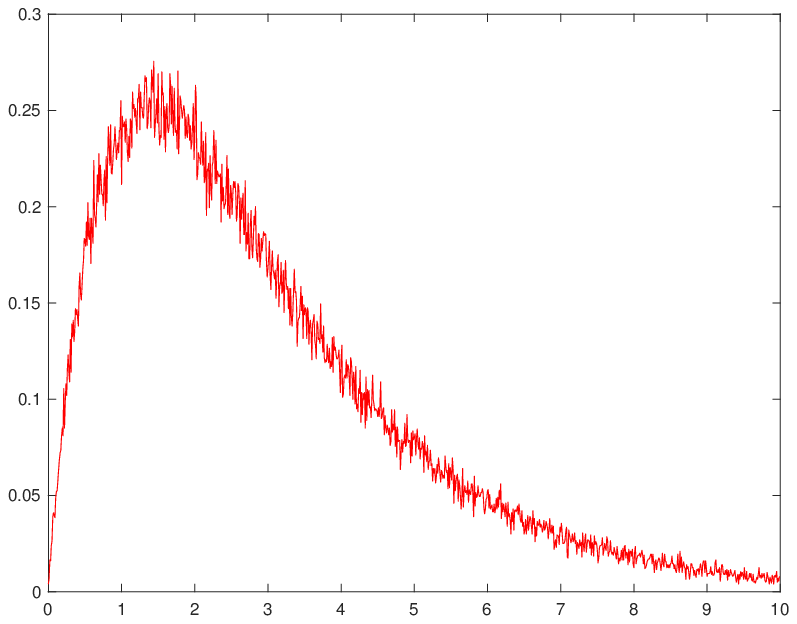}}
\caption{For $\theta = 1$ and $x = 2, $ estimated density of $\tau(x)$ (left) and  of $A(x)$ (right).}
\end{figure}

\subsection{Joint moments of $\tau (x)$ and $A(x).$ }
For $\lambda = (\lambda_1, \lambda _2 ) \in
\mathbb{R}_+ ^ 2 ,$ the joint Laplace transform of $(\tau(x), A(x))$ is
$$M_ \lambda (x) = M_ {(\lambda _1, \lambda _2)} (x) =E \left ( e^ { - \lambda _1 \tau (x) - \lambda _2 A(x) } \right ) $$
\begin{equation} \label{joint Laplace}
= E \left ( e^ { - \lambda _1 \tau (x) - \lambda _2 \int _0 ^ { \tau (x) }  X(t) dt } \right ) =
E \left ( e^ { - \int _0 ^ {\tau(x)} (\lambda _1 + \lambda _2 X(t) ) dt }\right ).
\end{equation}
By using \eqref{problemlaplace}
with $\lambda =1$ and $U(x)=
\lambda _1 + \lambda _2 x,$  we obtain that the function $M_
{(\lambda _1, \lambda _2 )} (x)$ solves the problem:
\begin{equation}\label{lab5}
\begin{cases}
\theta[M_{\lambda}(x - 1) - M_{\lambda}(x)] = (\lambda_{1}+\lambda_{2}x) M_{\lambda}(x), \;x>0; \\
M_{\lambda}(y) = 1, \; y \leq 0; \
\lim_{\mathbf{x} \to +\infty} M_{\lambda}(x) = 0.
\end{cases}
\end{equation}
The joint moments $E[\tau(x)^mA(x)^n]$, if they exist finite,
are given by:
\begin{equation}\label{lab6}
E[\tau(x)^mA(x)^n] = \left. (-1)^{m+n} \frac{\partial^m}{\partial \lambda_{1}^m}\frac{\partial^n}{\partial \lambda_{2}^n}
M_{\lambda}(x) \right|_{\lambda_{1} = \lambda_{2}= 0} .
\end{equation}
By taking the derivatives $\frac{\partial^2}{\partial\lambda_{1}\partial\lambda_{2}}|_{\lambda_{1} = \lambda_{2}= 0}$
in both sides of \eqref{lab5},
and equaling them, we obtain:
\begin{Theorem}\label{teo1}
$V(x) = E[\tau(x)A(x)]$ is the solution of the problem
\begin{equation}\label{lab7}
\begin{cases}
\theta[V(x - 1) - V(x)] = -xE[\tau(x)] - E[A(x)], \;x>0;\\
V(y) = 0, \; y \leq 0; \
\lim_{\mathbf{\theta} \to +\infty} V(x) = 0.
\end{cases}
\end{equation}
\end{Theorem}
\hfill  $\Box$
%\begin{Remark}
%The conditions for $y \leq 0$ and for $\theta \rightarrow +\infty$ follow from the facts that,
%if the starting point is  $y \leq 0$, then $\tau(y) = 0$ and $A(y) = 0$, while if the intensity $\theta$ of the Poisson process
%is $ \infty$, then the process $X(t) = x - N(t)$ reaches $0$ immediately.
%\end{Remark}
\par\noindent
Now we solve \eqref{lab7}, separately in the case when $x \in {\bf N}$, and when $x \not \in {\bf N}$.
\par\noindent
{\bf 1)} $x \in {\bf N}$:
from \eqref{primiduemomtaupoisson} and \eqref{mediaApoisson}, we obtain
$
-xE[\tau(x)] - E[A(x)] = -\frac{x^2}{\theta} - \frac{x(x+1)}{2\theta } = -\frac{x(3x+1)}{2\theta}.
$
Then, from the first equation of \eqref{lab7} we get
$ V(x) = V(x-1) + \frac{x(3x+1)}{2\theta^2}.$ It follows that
\begin{equation}\label{lab15}
V(x) = E[\tau(x)A(x)]
       = \sum_{i = 1}^x \frac{i(3i+1)}{2\theta^2} = \frac{x(x+1)^2}{2\theta^2},
\end{equation}
where  we used that
$\sum_{i=1}^n i = n(n+1)/2$ and $\sum_{i=1}^n i^2 = n(n+1)(2n+1)/6.$
The covariance of  $\tau(x)$ and $A(x)$ is
$$
Cov(\tau(x),A(x)) = E[A(x)\tau(x)] - E[A(x)]E[\tau(x)]
$$
\begin{equation}
= \frac{x(x+1)^2}{2\theta^2} - \frac{x}{\theta}\cdot\frac{x(x+1)}{2\theta} = \frac{x(x+1)}{2\theta^2}.
\end{equation}
The correlation coefficient  is obtained  by  \eqref{primiduemomtaupoisson} and
\eqref{secondomomApoisson}:
\begin{equation}
\rho(x) = \frac{Cov(\tau(x),A(x))}{\sqrt{Var(A(x))Var(\tau(x))}}
           = \sqrt{\frac{3(x+1)}{2(2x+1)}} .
\end{equation}
{\bf 2)} $x \not \in {\bf N}$:
from \eqref{primiduemomtaupoisson} and \eqref{mediaApoisson} we have
$
-xE[\tau(x)] - E[A(x)] = - \frac{([x]+1)(4x-[x])}{2\theta};
$
therefore, from the first  of \eqref{lab7} we obtain that
$
V(x) = V(x-1) + \frac{([x]+1)(4x-[x])}{2\theta^2}.
$
As in the previous case, it follows that
$V(x) = \sum_{n = 0}^{[x]} \frac{([x-n]+1)(4(x-n)-[x-n])}{2\theta^2} .$
Since $\forall k \in {\bf Z}$ and $ x \in \R,$  $[x + k] = k + [x],$
after calculation, we obtain:
\begin{equation}\label{lab16}
V(x) = \frac{\left([x]+1\right)\left([x]+2\right)\left(2x-[x]\right)}{2\theta^2} .
\end{equation}
Therefore,
$$Cov(\tau(x),A(x)) =
\frac{\left([x]+1\right)\left([x]+2\right)\left(2x-[x]\right)}{2\theta^2} -
\frac{\left([x]+1\right)}{\theta}\cdot\frac{(2x-[x])([x]+1)}{2\theta}  $$
\begin{equation}
= \frac{\left( [x] + 1 \right) \left( 2x - [x] \right)}{2\theta^2} .
\end{equation}
The correlation coefficient is:
$$\rho(x) = \frac{Cov(\tau(x),A(x))}{\sqrt{Var(A(x))Var(\tau(x))}} $$
\begin{equation}
= \frac{\frac{\left( [x] + 1 \right) \left( 2x - [x] \right)}{2\theta^2} }
{\sqrt{\frac{\left([x] + 1 \right)}{\theta^2}\cdot\frac{\left( [x] + 1 \right)}{12\theta^2}\left\lbrace 12x
\left( x - [x] \right) + 2[x] \left(2[x] + 1 \right) \right\rbrace}} =
\sqrt{\frac{3\left( 2x - [x] \right)^2}{12x\left( x - [x] \right) + 2x\left( 2[x] + 1 \right)}}.
\end{equation}
Notice that, both in the case when $x$ is an integer, and when it is not, $\rho(x)$ depends only on $x$ and,
$\rho(x) \rightarrow \sqrt{\frac{3}{4}} \simeq 0.8666 \ , $ as  $x\rightarrow +\infty .$\par\noindent

\begin{figure}[ht]
\centering
\includegraphics[width=%
0.9\columnwidth]{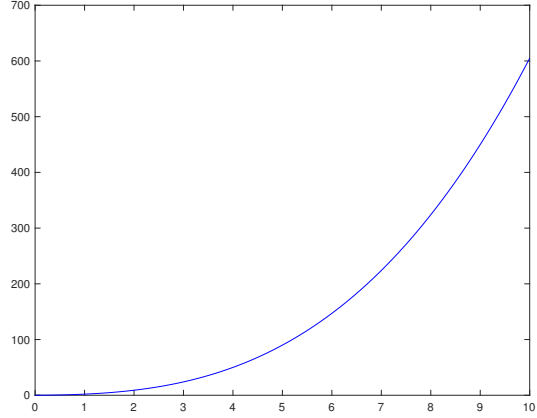}
\caption{Graph of $V(x)$ as a function of  $x\geq 0, $ for $\theta = 1$.
\label{fig1}}
\end{figure}

\subsection{Joint Laplace Transform of $(\tau(x),A(x))$}
Now, we find an explicit solution to the problem \eqref{lab5} for the joint Laplace transform of $\tau (x),  \ A(x), $ by
analyzing separately the case when $x \in {\bf N}$ and  $x \not
\in {\bf N}$. \par\noindent
From the first equation of \eqref{lab5}, we obtain $M_{\lambda}(x) =
\frac{\theta}{\lambda_{1}+\lambda_{2}x+\theta}M_{\lambda}(x-1).$
\par\noindent
{\bf 1)} $x \in {\bf N}:$ as easily seen, one has
\begin{equation}\label{lab17}
M_{\lambda}(x) =
\theta^x \cdot \prod_{k=1}^{x} \frac{1}{\lambda_{1}+\lambda_{2}k+\theta}.
\end{equation}
%We observe that the condition $\lim_{\mathbf{x} \to +\infty} M_{\lambda}(x) = 0$ is satisfied; actually, we have
%$$ (\lambda_{1}+\lambda_{2}x+\theta)(\lambda_{1}+\lambda_{2}(x-1)+\theta)\cdots(\lambda_{1}+\lambda_{2}+\theta) > \lambda_{2}x(x-1)(x-2)\cdots 2 \cdot 1$$
%So,
%\begin{equation}
%M_{\lambda}(x) < \frac{\theta^x}{\lambda_{2}x(x-1)(x-2)\cdots 2 \cdot 1}
%                         \thickapprox \frac{\theta^x}{x^x}\cdot \frac{1}{\lambda_{2}}
%= e^{x(\lg \theta - \lg x)}\cdot \frac{1}{\lambda_{2}} \longrightarrow 0, \quad
%                        for \quad x \rightarrow +\infty .
%\end{equation}
\noindent {\bf 2)} $x \not \in {\bf N}:$ if $x \in (0,1),$ then  $ M_{\lambda}(x) =
\theta / [\lambda_{1}+\lambda_{2}x+\theta],$ because $M_ \lambda (x-1)=1;$ if $x \in (1,2),$
then   $ M_{\lambda}(x) =
\theta^2 / [(\lambda_{1}+\lambda_{2}x+\theta)(\lambda_{1}+\lambda_{2}(x-1)+\theta)].$
\par\noindent
Thus, iterating  the procedure, one gets:
\begin{equation}\label{lab18}
M_{\lambda}(x) = \theta ^{ [x]+1} \prod _{k=0} ^{[x]} \frac 1 {\lambda _1 + \lambda _2 (x-k) + \theta } .
\end{equation}

\begin{figure}[ht]
\centering
\includegraphics[width=%
0.9\columnwidth]{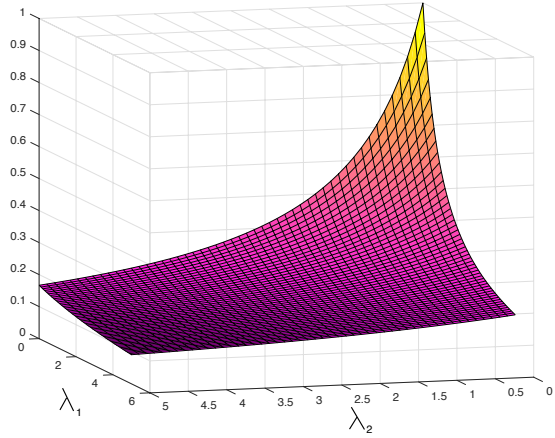}
\caption{Joint Laplace Transform of $(\tau(x),A(x)),$  as a function of $(\lambda _1, \lambda _2),$
in the case of Poisson process, for $x=1.$}
\end{figure}

It is easy to check that $\left. \frac{\partial^2}{\partial
\lambda_{1} \partial \lambda_{2}}M_{\lambda}(x)
\right|_{\lambda_{1} = \lambda_{2}= 0} = V(x)  $ where $V(x)$ is
given by ($\ref{lab15}$) and  \eqref{lab16}. In the Figure 6, we
show, in the case of Poisson process, the graph of the joint Laplace transform of
$(\tau (x),A(x)),$ as a function of $(\lambda _1, \lambda _2 ),$ for $x=1.$

\subsection{Moments of any order}
The aim of this subsection is to find the quantity
$V_{m, n}(x) := E \left[ \tau(x)^m A(x)^n \right], $
with $m, n \in {\bf N}$, following the same arguments used in
\cite{abdelve:mcap17} in the corresponding case.
Taking the \textit{m}-th partial derivative with respect to $\lambda_1$ and the \textit{n}-th partial derivative with
respect to $\lambda_2$ in the first equation of \eqref{lab5} and evaluating at $\lambda_1 = \lambda_2 = 0$, we obtain
the following:
\begin{Theorem}
 $V_{m, n}(x)$ is the solution of
\begin{equation}\label{nmoments}
\begin{cases}
\theta \left[ V_{m, n}(x - 1) - V_{m, n}(x) \right] = - nx V_{m, n - 1}(x) - m  V_{m - 1, n}(x) \\
V_{m, n}(y) = 0, \quad y \leq 0 ,  \
\lim_{\theta \to +\infty} V_{m, n}(x) = 0.
\end{cases}
\end{equation}
\end{Theorem}

\par \hfill  $\Box$
\begin{Remark}
For $m = n = 1$, we have $V_{1, 1}(x) = V(x) = E \left[ \tau(x) A(x) \right] $, $V_{0, 1}(x) = E \left[ A(x) \right]$
and $V_{1, 0}(x) = E \left[ \tau(x) \right] $, so \eqref{nmoments} becomes \eqref{lab7}.
\end{Remark}
By solving \eqref{nmoments} for $m = 2, n = 1$, one finds:
\begin{equation}
V_{2, 1}(x) =
\begin{cases}
\frac{x(x + 1)^2 (x + 2)}{2 \theta^3}, \qquad x \in {\bf N}\\\\
\frac{([x] + 1)([x] + 2)([x] + 3)(2x - [x])}{2 \theta^3}, \qquad x \not \in {\bf N}
\end{cases}
\end{equation}
that is, if $x \in {\bf N}$, $V_{2, 1}(x)$ is a polynomial of degree 4, otherwise it
has polynomial growth of degree $4.$
Similar polynomial expressions can be obtained for any $m$ and $n$, implying that the moments
$E \left[ \tau(x)^m A(x)^n \right]$ are finite, for all $m$ and $n$.
As far as the form of the solution $V_{m, n} (x)$ of \eqref{nmoments} is concerned,  proceeding by induction,
as done in \cite{abdelve:mcap17}, one gets:
\begin{Theorem}\label{teoMom}
For  integers $m, n \geq 0$, the solution of \eqref{nmoments} vanishes at zero, and it is a polynomial of degree
$m + 2n$ if
$x \in {\bf N}$, otherwise it has a polynomial growth of degree $m + 2n$.
\end{Theorem}
\par \hfill  $\Box$
\bigskip

When $x$ is an integer, by proceeding as done in \cite{abdelve:mcap17} for the analogous situation,
it is possible to obtain a compact closed  form of $V_{m,n} (x);$  in fact,
by using Theorem \ref{teoMom} for $x , m, n \in {\bf N}$, we obtain that there exist $a_1^{(m, n)}, \dots a_{m+2n}^{(m, n)} \in \R$ such that
$V_{m, n}(x) = \sum_{k = 1}^{m + 2n} a_k^{(m, n)} x^k ,$
moreover, there exist real numbers $a_k^{(m - 1, n)},$ for $k = 1, \dots, m+2n - 1$ and
$a_k^{(m, n - 1)}$ for $k = 1, \dots, m+2n - 2$ such that:
\begin{equation}
V_{m - 1, n}(x) = \sum_{k = 1}^{m+2n-1} a_k^{(m - 1, n)}x^k
 = a_1^{(m - 1, n)}x + \sum_{k = 1}^{m + 2n - 2} a_{k+1}^{(m - 1, n)}x^{k + 1},
\end{equation}
and
\begin{equation}
V_{m, n - 1}(x) = \sum_{k = 1}^{m + 2n - 2} a_k^{(m, n - 1)} x^k.
\end{equation}
Now, we introduce
three matrices $A^{(m, n)}$, $B^{(m-1, n)}$ and $C^{(m, n-1)} \in \R^{(m + 2n) \times (m + 2n)}$, such that
\begin{equation}
A_{i, j}^{(m,n)} =
\begin{cases}
\binom{j}{i-1} (-1)^{j - i +1}\theta, \qquad j \geq i = 1, \dots, m+2n\\\\
0, \qquad \qquad \qquad \qquad {\rm otherwise},
\end{cases}
\end{equation}
 \begin{equation}
B^{(m-1,n)} =
\begin{pmatrix}
0 & \dots & 0\\
\vdots & \ddots & \vdots\\
0 & \dots & -mI_{m + 2n - 1}
\end{pmatrix} \ \
{\rm and } \ \
C^{(m, n - 1)} =
\begin{pmatrix}
0 & \dots & 0\\
0 & \dots & 0\\
\vdots & \ddots & \vdots\\
0 & \dots & -nI_{m + 2n - 2}
\end{pmatrix}
\end{equation}
where $I_k \in \R^{k \times k}$ is the identity matrix. Finally, we denote by $\underline{a}_{m, n} \in \R^{m + 2n} $
the vector of the coefficients of $V_{m, n} (x)$, i.e.
$\underline{a}_{m, n} = \left( a_1^{(m, n)}, \dots, a_{m+2n}^{(m, n)} \right)^T .$
\par\noindent
By using \eqref{nmoments}, and proceeding as in \cite{abdelve:mcap17}, one gets
the following matrix equation:
\begin{equation}
A^{(m, n)} \underline{a}_{m, n} = B^{(m-1, n)}
\begin{pmatrix}
0\\
\vdots \\
\underline{a}_{m-1, n}
\end{pmatrix}
+ C^{(m, n-1)}
\begin{pmatrix}
0\\
0\\
\vdots \\
\underline{a}_{m, n-1}
\end{pmatrix}
\end{equation}
Being $A^{(m, n)}$ invertible, it results
\begin{equation}\label{amn}
\underline{a}_{m, n} = \left( A^{(m, n)} \right)^{-1} B^{(m-1, n)}
\begin{pmatrix}
0\\
\vdots \\
\underline{a}_{m-1, n}
\end{pmatrix}
+ \left( A^{(m, n)} \right)^{-1}  C^{(m, n-1)}
\begin{pmatrix}
0\\
0\\
\vdots \\
\underline{a}_{m, n-1}
\end{pmatrix}
\end{equation}

Equation \eqref{amn} provides a recursive formula to find the coefficient  $a_k ^{(m,n)},$ for $x \in {\bf N }.$
Thanks to the fact that the involved matrices
are triangular, for  $m$ and $n$ not too large, this formula
represents a faster way to obtain the coefficients of the polynomial $V_{m,n}(x)$ than solving directly the equation \eqref{nmoments}.
\subsection{Expected value of the time average till the FPT}
We fill find a closed form for
$$E \left ( \frac {A(x)} { \tau (x)}  \right )= E \left ( \frac 1  { \tau (x)} \int _0 ^{ \tau (x)} X(t) dt  \right ),$$
that is the expected value of the time average of Poisson process $X(t)= x - N_t $ till its FPT below zero. \par\noindent
If $M_ \lambda (x)= M_ {\lambda _1, \lambda _2} (x)$ is the joint Laplace transform of $(\tau(x), A(x)) ,$
defined by  \eqref{joint Laplace},
we note that
$$E \left ( A(x) e^{- \lambda _1 \tau (x)} \right ) = - \lim _{ \lambda _2 \rightarrow 0^+ } \frac \partial { \partial \lambda _2 }
M_ {\lambda _1, \lambda _2} (x),$$
and
\begin{equation} \label{A/tau}
E \left ( \frac {A(x)} { \tau (x)}  \right ) = \int _0 ^{ + \infty } E \left ( A(x) e^{- \lambda _1 \tau (x)} \right ) d \lambda _1 \ .
\end{equation}
Thus, the calculation of $E(A(x)/  \tau(x) )$ is reduced to integrate
$- \lim _{ \lambda _2 \rightarrow 0^+ } \frac \partial { \partial \lambda _2 } M_ {\lambda _1, \lambda _2} (x)$
with respect to $\lambda _1 .$ \par\noindent
By taking the partial derivative with respect to $\lambda _2$ in \eqref{lab17} and  \eqref{lab18}, and evaluating it at  $\lambda _2 =0,$
we easily obtain:
\begin{equation}
E( A(x) e^ {- \lambda _1 \tau (x)} ) =
\begin{cases}
\frac {x(x+1)} {2 }  \frac {\theta ^ x} {(\lambda _1 + \theta ) ^{x+1} }  & {\rm \ if} \ x \in {\bf N}, \\ $$
\left ([x]+1 \right ) \left (x - \frac {[x] } 2 \right )  \frac {\theta ^ {[x]+1} } {(\lambda _1 + \theta ) ^{[x]+2} }
& {\rm \ if} \ x \notin {\bf N}.
\end{cases}
\end{equation}
In the calculations, we have used the formula $(\prod _{i=1}^n f_i )' = (\prod _{i=1}^n f_i) ( \sum _ {i=1}^n f_i' / f_i ),$ giving the derivative of the product of $n$ functions. \par\noindent
Finally, by calculating the integral, we get:
\begin{equation}
E \left ( \frac {A(x) } { \tau (x)} \right ) =
 \int _ 0 ^ \infty E( A(x) e^ {- \lambda _1 \tau (x)} ) d \lambda _1 =
\begin{cases}
\frac { x+1} 2 & {\rm if } \ x \in {\bf N} \\
x - \frac {[ x]} 2 & {\rm if } \ x \notin {\bf N} .
\end{cases}
\end{equation}
Notice that, for any $x >0, \ E[A(x) / \tau  (x)]$ turns out to be $\ge x/2.$

\section{The case of Poisson Process with drift} For $x > 0$,
let us consider the process
\begin{equation}\label{PoissonDrift}
X(t) = x - \mu t - N_t
\end{equation}
where $\mu > 0 .$ As before, we denote by $\tau (x)$ the FPT below zero and by $A(x)= \int _0 ^ { \tau (x)} X(t) dt $ the FPA.
The infinitesimal generator is
\begin{equation} \label{generatorpoissondrift}
Lf(x) = - \mu \frac{\partial f}{\partial x} + \theta[f(x - 1) - f(x)], \  f \in C^{ 1}_b(\R) .
\end{equation}

\begin{Remark}\label{remPoisson}
Set
$Y(t) = x - N_t ;$ the FPT
$\tau_ Y(x) = \inf \{ t > 0 : Y_t \leq 0 \}$ and the FPA $A_Y(x) = \int_0^{\tau_Y(x)} Y(t)  dt$ have already
studied in the previous section. Since
$X(t)  \leq Y(t) ,$
one gets
$\tau(x) \leq \tau_Y(x),$ and also
$A(x) = \int_0^{\tau(x)} X(t) \, dt \leq \int_0^{\tau_Y(x)} Y(t) \, dt = A_Y(x). $
\end{Remark}
Let $U$ be a functional of the process $X(t);$ as before, for $\lambda > 0$ we denote by
$ M_{\lambda}(x) = \E\left[e^{-\lambda\int_{0}^{\tau(x)}U(X_s) \,ds}\right] $
the Laplace transform of the integral $\int_{0}^{\tau(x)}U(X_s) \,ds$.
Then, from
\eqref{problemlaplace}) of Theorem \ref{laplacetransform}, with the generator $L$ given by \eqref{generatorpoissondrift},
it follows that $M_{\lambda}(x)$ satisfies the problem with outer conditions:
\begin{equation}
\begin{cases}\label{lab35}
- \mu M'_{\lambda}(x) + \theta[M_{\lambda}(x - 1) - M_{\lambda}(x)] = \lambda U(x) M_{\lambda}(x), \;x > 0;\\
M_{\lambda}(y) = 1, \;y \leq 0; \
\lim_{x \to +\infty} M_{\lambda}(x) = 0.
\end{cases}
\end{equation}

\subsection{Laplace transform of $\tau(x)$} Taking $U(x) = 1$ in
\eqref{lab35}, we get that the Laplace transform of $\tau (x)$
satisfies:
\begin{equation}\label{lab36}
\begin{cases}
-\mu M'_{\lambda}(x) + \theta[M_{\lambda}(x - 1) - M_{\lambda}(x)] = \lambda M_{\lambda}(x), \; x > 0;\\
M_{\lambda}(y) = 1, \; y \leq 0; \
\lim_{\mathbf{x} \to +\infty} M_{\lambda}(x) = 0.
\end{cases}
\end{equation}
\noindent
{\bf 1)} if $x \in (0, 1],$ since
$M_{\lambda}(x - 1) = 1,$
the first equation of \eqref{lab36} becomes
$ \mu M'_{\lambda}(x) + \left( \theta + \lambda \right) M_{\lambda}(x) = \theta.$
By solving and taking into account the conditions of \eqref{lab36},  we find that
the Laplace transform of $\tau(x)$ is
$M_{\lambda}(x) = \frac{\lambda}{\theta + \lambda}e^{- \frac{\theta + \lambda}{\mu}x} + \frac{\theta}{\theta + \lambda} .$
Therefore,
$\E \left[ \tau(x) \right] = \left. - \frac{\partial}{\partial \lambda}M_{\lambda}(x) \right|_{\lambda = 0} = \frac{1}{\theta} - \frac{1}{\theta} e^{-\frac{\theta}{\mu}x} .$
\par\noindent
{\bf 2)} if $ x \in (1, 2],$
one has
$M_{\lambda}(x - 1) =\frac{\lambda}{\theta + \lambda}e^{- \frac{\theta + \lambda}{\mu}(x-1)} + \frac{\theta}{\theta + \lambda} .  $
Then, the first equation of \eqref{lab36} becomes
$ \mu M'_{\lambda}(x) + \left( \theta + \lambda \right) M_{\lambda}(x) =  \frac{\lambda \theta}{\theta + \lambda}e^{- \frac{\theta + \lambda}{\mu}(x-1)} + \frac{\theta^2}{\theta + \lambda}.$
By solving  and taking into account the conditions,  we find that
the Laplace transform of $\tau(x)$ is
$$M_{\lambda}(x) =  \frac{\lambda(\lambda + 2 \theta)}{(\theta + \lambda)^2}
e^{- \frac{\theta + \lambda}{\mu}x} +  \frac{\theta \lambda }{\mu (\theta + \lambda)} x e^{-\frac{\theta + \lambda}{\mu} (x - 1)} +
\frac{\theta^2}{(\theta + \lambda)^2} .$$
Therefore:
$$E \left[ \tau(x) \right] = \left. - \frac{\partial}{\partial \lambda}M_{\lambda}(x) \right|_{\lambda = 0} = \frac{2}{\theta} - \frac{2}{\theta}e^{-\frac{\theta}{\mu}x} + \frac{x}{\mu}e^{-\frac{\theta}{\mu}(x - 1)}.$$
\par\noindent
{\bf 3)} if $x \in (2, 3], $
one has $x - 1 \in (1, 2] $ and so
$M_{\lambda}(x-1) = \frac{\lambda(\lambda + 2 \theta)}{(\theta + \lambda)^2}e^{- \frac{\theta + \lambda}{\mu}(x - 1)} +  \frac{\theta \lambda }{\mu (\theta + \lambda)} (x - 1) e^{-\frac{\theta + \lambda}{\mu} (x - 2)} + \frac{\theta^2}{(\theta + \lambda)^2}.$
Then, the first equation of \eqref{lab36} becomes:
$$  \mu M'_{\lambda}(x) + \left( \theta + \lambda \right)M_{\lambda}(x) =$$
\begin{equation}\label{lab39}
\frac{\lambda \theta(\lambda + 2 \theta)}{(\theta + \lambda)^2}e^{- \frac{\theta + \lambda}{\mu}(x - 1)} +
\frac{\theta^2 \lambda }{\mu (\theta + \lambda)} (x - 1) e^{-\frac{\theta + \lambda}{\mu} (x - 2)} +
\frac{\theta^3}{(\theta + \lambda)^2}.
\end{equation}
By solving and taking into account the conditions of \eqref{lab36},  we find that
the Laplace transform of $\tau(x)$ is
$$M_{\lambda}(x) = \left[ 1 -  \frac{\theta^3}{(\theta + \lambda)^3} \right] e^{-\frac{\theta + \lambda}{\mu} x } + \frac{\theta \lambda (\lambda + 2 \theta)}{\mu^2 (\theta + \lambda)^2}x e^{-\frac{\lambda + \theta}{\mu}(x - 1)} $$
$$+ \frac{\theta^2 \lambda}{2 \mu^2(\theta + \lambda)}(x^2 - 2x)e^{-\frac{\lambda + \theta}{\mu}(x - 2)} +
\frac{\theta^3}{(\theta + \lambda)^3}.$$
In conclusions, for $x \in [0,3):$
\begin{equation} \label{LTtaupoissondrift}
M_ \lambda (x) =
\begin{cases}
\frac{\lambda}{\theta + \lambda}e^{- \frac{\theta + \lambda}{\mu}x} + \frac{\theta}{\theta + \lambda}, \ x \in [0,1) \\
\frac{\lambda(\lambda + 2 \theta)}{(\theta + \lambda)^2}e^{- \frac{\theta + \lambda}{\mu}x} +  \frac{\theta \lambda }{\mu (\theta + \lambda)} x e^{-\frac{\theta + \lambda}{\mu} (x - 1)} + \frac{\theta^2}{(\theta + \lambda)^2}, \ x \in [1,2) \\
\left[ 1 -  \frac{\theta^3}{(\theta + \lambda)^3} \right] e^{-\frac{\theta + \lambda}{\mu} x } + \frac{\theta \lambda (\lambda + 2 \theta)}{\mu^2 (\theta + \lambda)^2}x e^{-\frac{\lambda + \theta}{\mu}(x - 1)} \\
+ \frac{\theta^2 \lambda}{2 \mu^2(\theta + \lambda)}(x^2 - 2x)e^{-\frac{\lambda + \theta}{\mu}(x - 2)} +
\frac{\theta^3}{(\theta + \lambda)^3}, \ x \in [2,3)
\end{cases}
\end{equation}
Notice that, letting $\mu$ go to 0,  \eqref{LTtaupoissondrift} becomes \eqref{LTtaupoisson}, which holds in the case of Poisson process.

\begin{figure}[ht]
\centering
\includegraphics[width=%
0.9\columnwidth]{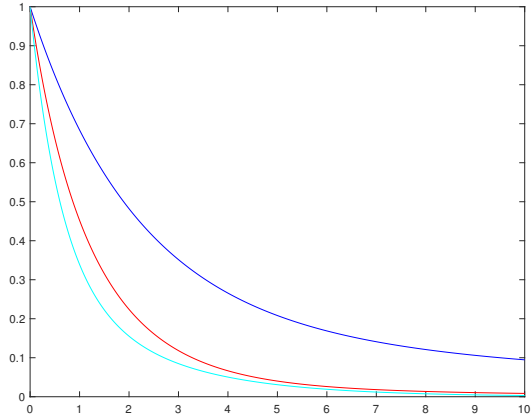}
\caption{Laplace Transform of $\tau(x),$ as a function of $ \lambda > 0 ,$ for $\mu = 1, \theta = 1;$  from top to bottom: $x = 0.5$ (blue), $x = 1.8$ (red) and $x = 2.5$ (cyan).}\label{fig:LaplaceTrasfDrift}
\end{figure}
In the  \figurename~\ref{fig:LaplaceTrasfDrift}, we report the shapes of the Laplace transform of
$\tau(x)$ for $\mu = 1$, $\theta = 1$ and various values of the starting point $x \in [0,3].$
\par\noindent
Owing to the complexity of the form, unlike the case of Poisson process, we are not able to find, neither a closed
form for the Laplace transform of $\tau(x)$ for any $x >0,$  nor a closed form for $E \left[ \tau(x) \right].$  However,
the upper and lower bounds hold:
\begin{equation} \label{inttau}
\frac{x}{\mu + \theta} \leq E \left[ \tau(x) \right] < \frac{x + 1}{\mu + \theta}
\end{equation}
In fact, at time $\tau(x)$, we have
$- 1 < x - \mu \tau(x) - N_{\tau(x)} \leq 0 $
for definition of $\tau(x),$ and because the size of the Poisson jumps is 1.
Taking expectation, we get
$- 1 < x - \mu E \left[\tau(x) \right] - E \left[N_{\tau(x)} \right]\leq 0.$
Since $E[N_ { \tau (x)}] = \theta E( \tau (x)),$
one obtains
$- 1 < x - \mu E \left[ \tau(x) \right] - \theta E \left[ \tau(x) \right] \leq 0 . $
Thus, from the left inequality, we obtain
$ E \left[ \tau(x) \right] < \frac{x + 1}{\mu + \theta} ;$
from the right inequality, we get
$ E \left[ \tau(x) \right] \geq \frac{x}{\mu + \theta}  .$
\begin{Remark}
$\tau (x) \le \tau _Y(x),$ where
$\tau_Y(x) = \inf \{ t > 0 : Y(t) \leq 0 \} , \ Y(t)=x - N_t .$
\end{Remark}

\subsection{Laplace transform of $A(x)$}
\begin{figure}[ht]
\centering
\includegraphics[width=%
0.9\columnwidth]{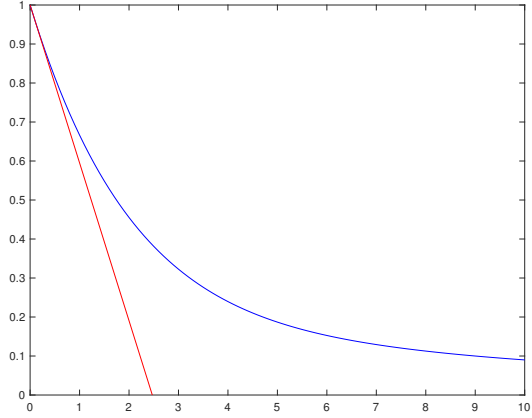}
\caption{Laplace transform of  $A(x),$ as a function of $\lambda > 0,$ obtained numerically by the Euler method with $x=1$, $\theta = 1$ and $\mu = 1$ and a discretization step $h = 0.1$ \ .}\label{fig:EuleroAxDrift}
\end{figure}

Taking $U(x) = x$ in
\eqref{lab35} , we get that the Laplace transform of $A (x)$
satisfies:
\begin{equation}
\begin{cases} \label{lab41}
- \mu M'_{\lambda}(x) + \theta[M_{\lambda}(x - 1) - M_{\lambda}(x)] = \lambda x  M_{\lambda}(x), \;x>0;\\
M_{\lambda}(y) = 1, \; y \leq 0; \
\lim_{x \to +\infty} M_{\lambda}(x) = 0.
\end{cases}
\end{equation}
As before, we look for a solution, taking $x >0$ in successive intervals.
We start with $x \in (0, 1]$:
it results $ x - 1 \leq 0, $ and so
$ M_{\lambda}(x - 1) = 1.$
Then, the first equation of \eqref{lab41} becomes
$ \mu M'_{\lambda}(x) + \left( \theta + \lambda x \right) M_{\lambda}(x) = \theta ,$ whose
solution is:
$$M_{\lambda}(x) = exp\biggl\{- \left( \frac{\theta}{\mu}x + \frac{\lambda}{\mu}\frac{x^2}{2} \right)\biggr\} \left[ c + \frac{\theta}{\mu}
\frac{1}{\sqrt{\lambda}} e^{-\frac{\theta^2}{2 \lambda \mu}} \int _0 ^x  e^{\frac{t^2}{2\mu}} \ dt \right ] ,
$$
where c is a constant. Since the last integral cannot be written in terms of elementary functions, we cannot iterate this procedure, to find the
solution of \eqref{lab41} for every $x >0.$
However, it is possible to integrate numerically \eqref{lab41} by the Euler method; for $x=1,$
the shape of $M_{\lambda} (x),$ so obtained, as a function of $\lambda >0 ,$  is reported in \figurename~\ref{fig:EuleroAxDrift}.

We have  found $E[A(x)]$ as the value, changed in sign, of the slope of the tangent line at zero of the Laplace transform of $A(x).$
Then, we have estimated $E \left[ A(x) \right] $ also by Monte Carlo simulation.
In \figurename~\ref{fig:EuleroAxDrift}, the slope of the tangent line at zero of the Laplace transform is
$m = -0.404$ and so, for $x=1$, $\theta = 1$ and $\mu = 1$, we obtain the estimate $\widehat { E[A(x)]} = 0.404$, while
Monte Carlo simulation provides the  value $0.399$.

\begin{Remark}
Of course, the FPA of Poisson process with drift is less than the one regarding  Poisson process. \par\noindent
Notice that, neither the joint moments of $\tau(x)$ and $ A (x) ,$  nor their joint Laplace transform,
can be obtained in closed form, because the corresponding PDDEs cannot be solved explicitly.
\end{Remark}
\section{The case of Drifted Brownian motion with Poisson
jumps} For $x > 0, \ \mu >0,$ and $\sigma >0,$  let be
$X(t) = x - \mu t + \sigma B_t - N_t .$
The infinitesimal generator is

%$$ Lg(x) = - \mu \frac{\partial}{\partial x} + \frac{1}{2} \frac{\partial^2}{•\partial x^2}}+
%\theta[g(x - 1) - g(x)]}, \; g \in C^{ 2}_b(\R) $$
\begin{equation} \label{generatordriftedBMpoisson}
Lf(x) = - \mu \frac {\partial f}{\partial x} + \frac{\sigma ^2 }{2} \frac {\partial ^2 f} {\partial x^2 }+
\theta[f(x - 1) - f(x)], \ f \in C^{ 2}_b(\R)
\end{equation}

%\begin{figure}[ht]
%\centering
%\includegraphics[width=%
%0.9\columnwidth]{GraficoPoissonDriftBt.eps}
%\caption{A sample path of the process $X(t) = x - \mu t + B_t-  N(t)$ with starting point $x=10$, $\theta = 1$ and $\mu = 1$.}
%\end{figure}

As before, let $U$ be a functional of the process $X$ and let
$ M_{\lambda}(x) = E\left[e^{-\lambda\int_{0}^{\tau(x)}U(X_s) \,ds}\right] $
the Laplace transform of the integral $\int_{0}^{\tau(x)}U(X_s)
\,ds$.
Then, from  \eqref{problemlaplace} of Theorem \ref{laplacetransform}, with generator $L$ given by \eqref{generatordriftedBMpoisson},
one gets that $M_{\lambda}(x)$ satisfies
the problem with outer conditions:
\begin{equation}\label{lab57}
\begin{cases}
\frac{\sigma ^2 }{2}M''_{\lambda}(x) - \mu M'_{\lambda}(x) + \theta[M_{\lambda}(x - 1) - M_{\lambda}(x)] = \lambda U(x) M_{\lambda}(x), \;x > 0;\\
M_{\lambda}(y) = 1, \;y \leq 0; \
\lim_{x \to +\infty} M_{\lambda}(x) = 0.
\end{cases}
\end{equation}
\subsection{Laplace transform of $\tau(x)$}
Taking $U(x) = 1$ in (\ref{lab57}), we
get that the Laplace transform of $ \tau (x)$ is the solution of the PDDE:
\begin{equation}\label{lab58}
\begin{cases}
\frac{\sigma ^2 }{2}M''_{\lambda}(x) - \mu M'_{\lambda}(x) + \theta[M_{\lambda}(x - 1) - M_{\lambda}(x)] = \lambda M_{\lambda}(x), \;x > 0;\\
M_{\lambda}(y) = 1, \;y \leq 0; \
\lim_{x \to +\infty} M_{\lambda}(x) = 0;
\end{cases}
\end{equation}
{\bf 1)} if $x \in (0, 1],$ then $M_{\lambda}(x - 1) = 1 ;$
the first equation of (\ref{lab58}) becomes
$\frac{\sigma ^2}{2}M''_{\lambda}(x) - \mu M'_{\lambda}(x) + \theta[1 - M_{\lambda}(x)] = \lambda M_{\lambda}(x). $
By solving and taking into account the conditions of \eqref{lab58}, we find that
the Laplace transform of
$\tau(x)$ is
\begin{equation} \label{LTtaudriftedBM+poissonfirstinterval}
M_{\lambda}(x) = \frac{\lambda}{\lambda +  \theta}e^{\left( \mu / \sigma ^2 - (1 / \sigma ^2) \sqrt{\mu^2 + 2(\lambda + \theta )} \right) x } +  \frac{\theta}{\lambda + \theta}.
\end{equation}
\noindent {\bf 2)} if $x \in (1, 2],$
one has
$M_{\lambda}(x-1) = \frac{\lambda}{\lambda +  \theta}e^{\left( \mu / \sigma ^2 - (1 / \sigma ^2) \sqrt{\mu^2 + 2(\lambda + \theta )} \right) x } +  \frac{\theta}{\lambda + \theta}.$
Then, the first equation of (\ref{lab58}) becomes:
$$
\frac{\sigma ^2 }{2}M''_{\lambda}(x) - \mu M'_{\lambda}(x) - (\lambda + \theta)M_{\lambda}(x) = - \frac{\lambda \theta}{ \sigma ^2 (\lambda + \theta)}e^{\left(
\mu / \sigma ^2 - (1 / \sigma ^2) \sqrt{\mu^2 + 2 \sigma ^2 (\lambda + \theta )} \right) (x - 1) }   \frac{\theta^2}{ \sigma ^2 (\lambda + \theta)}.
$$
By solving and taking into account the
conditions, we find that
the Laplace transform of
$\tau(x)$ is
$$M_{\lambda}(x) =  \frac{\lambda(\lambda + 2 \theta)}{\sigma ^2 (\lambda + \theta)} e^{\left( \mu / \sigma ^2 - (1/ \sigma ^2) \sqrt{\mu^2 + 2 \sigma ^2 (\lambda + \theta )} \right) x } + $$
\begin{equation} \label{LTtaudriftedBM+poissonsecondinterval}
\frac{\sigma ^4 \lambda \theta}{ (\lambda + \theta) \sqrt{\mu^2 + 2 \sigma ^2 (\lambda + \theta)}}x e^{\left( \mu / \sigma ^2 - (1 / \sigma ^2) \sqrt{\mu^2 +
2 \sigma ^2 (\lambda + \theta )} \right) (x - 1)} + \left( \frac{\theta}{\lambda + \theta} \right)^2 .
\end{equation}
{\bf 3)} if $x \in (2, 3],$
one has
$M_{\lambda}(x-1) =  \frac{\lambda(\lambda + 2 \theta)}{\sigma ^2 (\lambda + \theta)} e^{\left( \mu / \sigma ^2 - 1/ \sigma ^2 \sqrt{\mu^2 + 2 \sigma ^2 (\lambda + \theta )} \right) (x-1) } + $ \par\noindent
$+\frac{\sigma ^4 \lambda \theta}{(\lambda + \theta) \sqrt{\mu^2 + 2 \sigma ^2 (\lambda + \theta)}}(x-2) e^{\left( \mu / \sigma ^2 - (1 / \sigma ^2) \sqrt{\mu^2 +
2 \sigma ^2 (\lambda + \theta )} \right) (x - 2)} + \left( \frac{\theta}{\lambda + \theta} \right)^2 .$
Then, the first equation of (\ref{lab58}) becomes:
$$\frac{\sigma ^2 }{2}M''_{\lambda}(x) - \mu M'_{\lambda}(x) -(\lambda + \theta)M_{\lambda}(x) =  -\frac{\lambda \theta(\lambda + 2 \theta)}{ \sigma ^4 \lambda + \theta} e^{\left( \mu / \sigma ^2 - (1/ \sigma ^2) \sqrt{\mu^2 + 2 \sigma ^2 (\lambda + \theta )} \right) (x-1) } + $$
$$-\frac{\lambda \theta^2}{\sigma ^2 (\lambda + \theta) \sqrt{\mu^2 + 2 \sigma ^2 (\lambda + \theta)}}(x-1)
e^{\left( \mu / \sigma ^2 - (1/ \sigma ^2) \sqrt{\mu^2 + 2 \sigma ^2 (\lambda + \theta )} \right) (x-2) } -  \frac{\theta^3}{\sigma ^2 (\lambda + \theta)^2}.$$
By solving and taking into account the
conditions, we find that
the Laplace transform of
$\tau(x)$ is
$$M_{\lambda}(x) = (1 - (\theta / (\lambda + \theta ))^3 ) e^{\left( \mu / \sigma ^2 - (1/ \sigma ^2) \sqrt{\mu^2 + 2 \sigma ^2 (\lambda + \theta )} \right) x }
$$
$$+
\frac {\theta^2 \lambda (\lambda + 2 \theta) } {\sigma ^4 (\theta + \lambda)^2\sqrt{\mu^2 + 2(\theta + \lambda)}}
x  e^{\left( \mu / \sigma ^2 - (1/ \sigma ^2) \sqrt{\mu^2 + 2 \sigma ^2 (\lambda + \theta )} \right) (x-1) }
$$
$$ + \frac{\theta^2 \lambda }{2(\theta + \lambda)(\mu^2 + 2 \sigma ^2 (\theta + \lambda))} \left( x^2 + 2x ( \frac{\sigma ^2 }{2\sqrt{\mu^2 + 2(\lambda + \theta)}} - 1 ) e^{\left( \mu / \sigma ^2 - (1/ \sigma ^2) \sqrt{\mu^2 + 2 \sigma ^2 (\lambda + \theta )} \right) (x-2) } \right) $$
\begin{equation} \label{LTtaudriftedBM+poissonthirdinterval}
+ \left( \frac{\theta}{\theta + \lambda} \right) ^3 .
 \end{equation}

In line of principle, it is possible to iterate this procedure to obtain the Laplace transform of $\tau (x),$ also for $ x>3,$ but the calculations are very
complicated. Therefore, for $x >3$ it is more convenient to estimate numerically the
Laplace transform of $\tau(x),$ and $E \left[ \tau(x) \right],$
by Monte Carlo simulations.

\begin{figure}[ht]
\centering
\includegraphics[width=%
0.9\columnwidth]{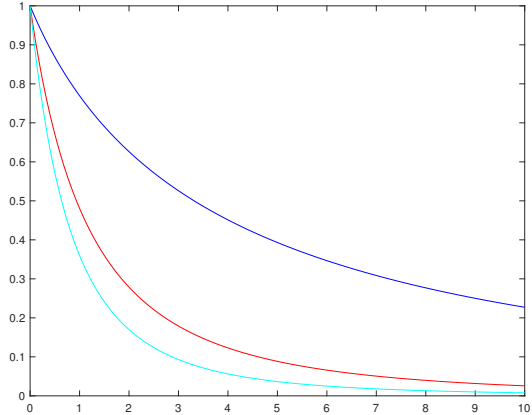}
\caption{Laplace Transform of $\tau(x),$ as a function of $ \lambda >0, $  for $\mu = 1, \theta = 1;$ from top to bottom:  $x = 0.5$ (blue), $x = 1.8$ (red) and $x = 2.5$ (cyan).}\label{fig:LaplaceTrasfDriftBt}
\end{figure}

In the  Figure 9, we  show the shapes of the Laplace transform of $\tau(x)$
for $\mu = 1 , \ \sigma =1,  \ \theta = 1$ and various values of the starting point $x \in [0,3],$ obtained by using \eqref{LTtaudriftedBM+poissonfirstinterval},
\eqref{LTtaudriftedBM+poissonsecondinterval} and \eqref{LTtaudriftedBM+poissonthirdinterval}. \par\noindent
By calculating the derivative of $M_{\lambda}(x)$ with respect to $\lambda ,$ at $ \lambda =0,$ one obtains $- E( \tau (x));$
for instance, for $x \in (0,1]$ one gets:
\begin{equation}
E( \tau (x))= \frac 1 \theta [ 1- e^{ \frac 1 { \sigma ^2} \left ( \mu - \sqrt {\mu ^2 + 2 \theta } \ \right )x} ] .
\end{equation}

\subsection{Laplace transform of $A(x)$} Taking $U(x) = x$ in
\eqref{lab57}, we get that the Laplace transform of $A(x)$ is the solution of the PDDE:
\begin{equation}\label{lab60}
\begin{cases}
\frac{1}{2}M''_{\lambda}(x) - \mu M'_{\lambda}(x) + \theta[M_{\lambda}(x - 1) - M_{\lambda}(x)] = \lambda x M_{\lambda}(x), \;x > 0;\\
M_{\lambda}(y) = 1, \;y \leq 0; \
\lim_{x \to +\infty} M_{\lambda}(x) = 0.
\end{cases}
\end{equation}
Unfortunately,
even for $x \in (0, 1]$, the first equation of \eqref{lab60} is a second order ODE with non-constant coefficients, so
its explicit solution cannot be found; neither the solution can be calculated numerically by the Euler method, since the initial condition  $M'_{\lambda}(0)$ is unknown. As far as the expectation of $A(x)$ is concerned,
one has to  estimate it by Monte Carlo simulation by using
$$A(x) = \int_0^{\tau(x)} X(t) \, dt = x \tau(x) - \mu \frac{\tau(x)^2}{2} + \int_0^{\tau(x)} B_t \, dt - \int_0^{\tau(x)} N_t \, dt. $$
For instance, for $x = 10$, $\theta = 1, \  \mu = 2, $ and $\sigma =1,$  Monte Carlo simulation provides the estimate
$\widehat {E[A(x)]} =  17.782.$
Notice, however, this
estimate is not the more accurate possible, but it has only
an illustrative purpose; in fact, we have not considered any procedure to avoid the overestimate of the FPT time $\tau (x)$ (see e.g. \cite{gir:mcap01}).

\section{Conclusions and final remarks}
In this paper,
we have carried on the study,  already undertook in \cite{abundo:mcap13} for general jump-diffusions, and in \cite{abdelve:mcap17}
for drifted Brownian motion,
of the first-passage area (FPA)  $A(x),$ swept out by a one-dimensional jump-diffusion process $X(t),$ starting from $x>0,$
till its first-passage time (FPT) $\tau(x)$
below zero.
Here, we have investigated the joint distribution of $\tau(x)$ and $A(x),$ in the special case when
$X(t )$ is a L$\acute{\text{e}}$vy process of the form $X(t)= x - \mu t + \sigma B_t - N_t,$ where $\mu \ge 0, \  \sigma \ge 0, \ B_t$
is standard Brownian motion, and $N_t$ is a homogeneous Poisson process with intensity $\theta,$ starting at zero.
We have established partial differential-difference equations (PDDE's) with outer conditions for the Laplace transform
of the random vector $(\tau(x), A(x)),$ and for
the joint moments  $E[\tau(x)^m A(x)^n]$ of the FPT and FPA. In the special case of Poisson process $( \mu = \sigma =0),$ that is $X(t)= x -N_t,$
we have presented an  algorithm  to find recursively
$E[\tau(x)^m A(x)^n],$ for any $m$ and $n;$
moreover, we have calculated the expected value of the time average of $X(t)$ till its FPT below zero.
In the other cases, we have not been always able to carry on  explicit calculations, whenever the corresponding PDDEs
equations cannot be solved in closed form. \par
We observe that,
for a given barrier $S$, one can extend the
results of this paper to the FPT and FPA through $S$ of a one-dimensional L$\acute{\text{e}}$vy process
$X(t)= x + \mu t + \sigma B_t + N_t \ (\mu , \sigma \ge 0),$ starting from $x <S;$ in this article, we have
considered the case of crossing zero,  for the sake of
simplicity. Really, for $x < S,$ one can consider
the FPT $\tau_S(x) = \inf \{ t > 0 : X(t) \ge S  \} $ and
the FPA $A_S(x) = \int_0^{\tau_S(x)} X(t) \, dt$ determined by
$X(t)$ till its FPT $\tau_S(x)$ through $S$; note that $A_S(x)$ is
improperly called the FPA of $X(t)$ through $S$, since the area of
the plane region determined by the trajectory of $X(t)$ and the
time axis in the first-passage period $[0, \tau_S(x)]$ is $
\int_0^{\tau_S(x)} |X(t)| \ dt $, which coincides with $A_S(x)$
only if $X(t)$ is non-negative in the entire interval $[0,
\tau_S(x)]$.
However, also in this case one obtains that the
Laplace transforms and moments of $\tau_S(x)$ and $A_S(x)$ are
solutions to certain PDDE's with outer conditions, which are
similar to those analyzed in the present paper.  In fact,
for $X(t)= x + \mu t + \sigma B_t + N_t ,$ with $x <S,$ one
has
$$ \tau _S(x) = \inf \{ t >0: x + \mu t + \sigma B_t + N_t \ge S  \} = \inf \{ t >0: -x - \mu t - \sigma B_t - N_t \le -S  \};$$
Since $- B_t \sim B_t,$ then $\tau _S(x)$ has the same distribution as
$$  \inf \{ t >0: S -x - \mu t + \sigma B_t - N_t \le 0  \} \equiv \widetilde \tau  (S-x),$$
where $\widetilde \tau (z) $ is the FPT of
L$\acute{\text{e}}$vy process $Z(t)= z - \mu t + \sigma B_t - N_t$
below zero, which was studied in this paper.
Moreover:
$$ A_S(x) = \int _ 0 ^ {\tau _S (x)} X(t) dt = S \tau (x) - \widetilde A(S-x) = S \widetilde \tau (S-x) - \widetilde A(S-x) ,$$
where
$$\widetilde A(z) = \int _0 ^ {\widetilde \tau  (z)} (z - \mu t + \sigma B_t - N_t ) dt $$
is the FPA below zero of drifted BM with Poisson jumps, already studied. This implies that the Laplace transforms and moments
of $\tau_S(x)$ and $A_S(x)$ are again solutions to suitable
PDDE's with outer conditions (see e.g. \cite{abundo:open13}). \par
In line of principle, all the results contained in this paper can be extended  to a jump-diffusion process $X(t)$
in which the drifted BM
is replaced by a one-dimensional, time homogeneous  diffusion, and also to a more
general L$\acute{\text{e}}$vy process, in which the jump sizes are not constant,
but they are independent random variables with assigned distribution,
implying that $N_t$ becomes a compound Poisson process; of course, explicit calculations become more complicated.

\end{document}